\documentclass[runningheads,a4paper]{llncs}
\usepackage{amssymb}
\usepackage{graphicx,color}
\usepackage[table,xcdraw]{xcolor}
\usepackage{wrapfig}
\usepackage{hyperref}
\usepackage{amssymb,amsmath,mathrsfs,dsfont,bbold,latexsym,enumerate,pstricks,epsf}%
\usepackage[font=footnotesize]{caption}
\usepackage[font=scriptsize]{subcaption}
\usepackage{url}
\usepackage{soul}
\urldef{\mailsa}\path|carole.le-guyader@insa-rouen.fr|
\urldef{\mailsb}\path|ai323@cam.ac.uk|
\urldef{\mailsc}\path|mjg40@cam.ac.uk|
\urldef{\mailsd}\path|{vc324, ai323, nd448, mjg40, cbs31, gbw1000}@cam.ac.uk|
\urldef{\mailse}\path|gbw1000@cam.ac.uk|
\newcommand{\keywords}[1]{\par\addvspace\baselineskip
\noindent\keywordname\enspace\ignorespaces#1}

\usepackage[colorinlistoftodos]{todonotes}
\makeatletter
\renewcommand\footnotesize{%
   \@setfontsize\footnotesize\@ixpt{11}%
   \abovedisplayskip 8\p@ \@plus2\p@ \@minus4\p@
   \abovedisplayshortskip \z@ \@plus\p@
   \belowdisplayshortskip 4\p@ \@plus2\p@ \@minus2\p@
   \def\@listi{\leftmargin\leftmargini
               \topsep 4\p@ \@plus2\p@ \@minus2\p@
               \parsep 2\p@ \@plus\p@ \@minus\p@
               \itemsep \parsep}%
   \belowdisplayskip \abovedisplayskip
}
\makeatother 

\makeatletter
\renewcommand\footnotesize{%
   \@setfontsize\footnotesize\@ixpt{7}%
   \abovedisplayskip 8\p@ \@plus2\p@ \@minus4\p@
   \abovedisplayshortskip \z@ \@plus\p@
   \belowdisplayshortskip 4\p@ \@plus2\p@ \@minus2\p@
   \def\@listi{\leftmargin\leftmargini
               \topsep 4\p@ \@plus2\p@ \@minus2\p@
               \parsep 2\p@ \@plus\p@ \@minus\p@
               \itemsep \parsep}%
   \belowdisplayskip \abovedisplayskip
}
\makeatother
\begin{document}

\mainmatter  

\title{Multi-tasking to Correct: Motion-Compensated MRI via Joint Reconstruction and Registration}

\titlerunning{}
%
%


\author{Veronica Corona\inst{1}\and Angelica I. Aviles-Rivero\inst{2}\and No\'emie Debroux\inst{1}\and Martin Graves\inst{3} \and Carole Le Guyader\inst{4}\and Carola-Bibiane Sch\"onlieb\inst{1} \and Guy Williams\inst{5} }

\authorrunning{V. Corona et al.}   
\titlerunning{ Motion-Compensated MRI via Joint Reconstruction and Registration}



\institute{
DAMTP, University of Cambridge,  United Kingdom\\
\and
DPMMS, University of Cambridge,  United Kingdom\\
\and
Department of Radiology, University of Cambridge, United Kingdom
\\
\and
LMI, Normandie Universit\'e, INSA de Rouen,  France\\
\mailsa \\
\and
Department of Clinical Neurosciences, University of Cambridge, United Kingdom
\mailsd
}

%
%

\maketitle
\hyphenation{va-ria-tional}
\hyphenation{quasicon-vexity}
\hyphenation{si-gnificant}
\hyphenation{analy-sis}
\hyphenation{facili-tate}
\hyphenation{di-sease}
\hyphenation{reg-is-tra-tion}
\hyphenation{diffe-rence}
\hyphenation{dissimila-rity}
\hyphenation{desi-gning}
\hyphenation{hyper-elastic}
\hyphenation{pro-per-ty}
\hyphenation{li-te-ra-tu-re}
\hyphenation{boun-da-ry}
\hyphenation{mi-ni-mi-zing}
\hyphenation{e-ve-ry}
\hyphenation{fea-tu-re}
\hyphenation{pe-na-li-zing}
\hyphenation{ne-ga-ti-ve}
\hyphenation{pro-blem}
\hyphenation{regu-larization}
\hyphenation{mo-del-ling}
\hyphenation{drop-ped}
\hyphenation{lo-ca-ting}
\hyphenation{ref-er-ence}
\hyphenation{min-i-mi-za-tion}
\hyphenation{em-pir-i-cal-ly}

\begin{abstract}
This work addresses a central topic in  Magnetic Resonance Imaging (MRI) which is the motion-correction problem in a joint reconstruction and registration framework. From a set of multiple MR acquisitions corrupted by motion, we aim at - \textit{jointly} - reconstructing a single motion-free corrected image and retrieving the physiological dynamics through the deformation maps. To this purpose,  
we propose a novel variational model. First, we introduce an $L^2$ fidelity term, which  intertwines  reconstruction and registration along with  the weighted total variation. Second, we introduce an additional regulariser which is based on the hyperelasticity principles to allow large and smooth deformations. We demonstrate through numerical results that this combination creates synergies in our complex variational approach resulting in higher quality reconstructions and a good estimate of the breathing dynamics. We also show that our joint model outperforms in terms of contrast, detail and blurring artefacts, a sequential approach.
\end{abstract}

\section{Introduction}\label{sec::introduction}
Magnetic Resonance Imaging (MRI) is a well-established modality that allows to capture details of almost any organ or anatomical structure of the human body. However, the lengthy period of time needed to acquire the necessary measurements to form an MR image continues to be a major contributor to image degradation, which compromises clinical interpretation~\cite{zaitsev2015motion}. This prolonged time makes MRI highly sensitive to motion mainly because the timescale of physiological motion -including involuntary motion- is generally shorter than the acquisition time~\cite{zaitsev2015motion}. This motion is manifested as undesirable artefacts including geometric distortions and blurring, which causes a significant degradation of the image quality and affects the clinical relevance for diagnosis. Therefore, the question of \textit{how to do motion correction in MRI?} is of great interest at both theoretical and clinical levels, and this is the problem addressed in this paper. 

Although motion reduction strategies such as faster imaging (e.g.~\cite{lustig2007sparse}) have shown potential results, \textit{we focus on post-acquisition correction of MRI sequences}, which involves explicit estimation of the motion. In this context, there have been different attempts for motion correction in MRI, whose predominant scheme is based on image registration after reconstructing the images (i.e. to estimate a mapping between a pair of images). A set of algorithmic approaches use rigid registration including~\cite{gupta2003fast,adluru2006model,johansson2018rigid}. 
However, the intrinsic nature of the problem yields those techniques to fail to capture complex transformations - e.g.  cardiac and breathing motion - which, in many clinical cases, are inevitable. To mitigate this limitation, different non-rigid image registration methods have been proposed in this domain - e.g.~\cite{ledesma2007motion,jansen2017evaluation,Burger2013,droske}. 

On the other side of the aforementioned methods that separately compute the image reconstruction and the motion estimation, there is another type of algorithmic approaches that computes those sub-tasks jointly. This philosophy was early explored in~\cite{tomasi1992shape}. In that work, authors proposed a factorisation method, which uses the singular-value decomposition for recovering shape and motion under orthography. Particularly in the medical domain, and following a variational perspective~\cite{burger2018variational}, different works have been reported including for SPECT imaging~\cite{mair2006estimation,schumacher2009combined}, PET~\cite{blume2010joint} and  MRI~\cite{aviles2018compressed}- to name a few. The works with a closer aim to ours are discussed next. 

Schumacher et al.~\cite{schumacher2009combined} presented an algorithmic approach that combines reconstruction and motion correction for SPECT imaging. The authors proposed a novel variational approach using a regulariser that penalises an offset of motion parameter to favour a mean location of the target object. However, the major limitation is that they only consider rigid motions.
In the same spirit, Fessler~\cite{fessler2009} 
proposed a generic joint reconstruction/registration framework with a model based on a penalised-likelihood functional using a weighted least square fidelity term along with a spatial and a motion regulariser. 
Authors in~\cite{blume2010joint} propose a joint model composed of a motion-aware likelihood function and a smoothing  term for a simultaneous image reconstruction and motion estimation for PET data. More recently, Odille et al~\cite{odille2016joint} proposed a joint model for MRI image reconstruction and motion estimation. This approach allows for an estimate of both intra and inter-image motion, meaning that, not only the misalignment problem is addressed but also it allows correcting for blurring/ghosting artefacts. 

Most recently and motivated by the deep learning (DL) performance breakthrough, a set of algorithmic approaches has been devoted to investigate the benefits of DL for image registration - e.g.~\cite{yang2017quicksilver,de2017end}. Although, certainly, those approaches deserve attention, their review goes beyond the scope of this paper.

\textbf{Contributions.} In this work, we follow the joint model philosophy, in which we seek to compute - simultaneously and jointly - the MR image reconstruction representing the true underlying anatomy and the registration tasks. Whilst the use of joint models (multi-tasking) has been widely explored for modalities such as PET and SPECT, the application for MRI is a largely untouched area. This is of great interest as 
MRIs highly sensitive to motion. With this motivation in mind we propose a variational model that allows for intra-scan motion correction in MRI. The significance of our approach is that by computing the reconstruction and registration tasks jointly one can exploit their strong correlation thus reducing error propagation and resulting in a significant motion correction. Whilst this is an important part of our solution, our main contributions are:

\begin{itemize}
    \item A novel theoretically well-motivated mathematical framework for motion correction in MRI. It relies on multi-tasking to address reconstruction and registration jointly, and guarantees the preservation of anatomical structures and topology by ensuring the positivity of the Jacobian determinant.
    \item Evaluation of our theory on a set of numerical and visual results using MRI datasets. We also show that our joint model leads to a better output than a sequential approach using the FAIR toolbox with similar characteristics. 
\end{itemize}

\section{Joint Model for Motion Correction in MRI}  
In this section, we  describe our approach for motion correction using MRI data. 
Firstly, we theoretically motivate our joint reconstruction and registration model and secondly, we formulate our joint variational model for motion correction. 
Finally, we describe an optimisation scheme to numerically solve our proposed joint model.

\subsection{Mathematical Model and Notation}
In an MRI setting, a target image $u$ representing a part of the patient body  is acquired in spatial-frequency space.  Let us denote by $x_i\in L^2(\mathbb{R}^2)$ the $i$-th  measurement acquired by an MRI scanner for $i=1,\cdots,T$, where $T$ is the number of acquisitions. The measured $x_i$ can be modelled as $x_i=\mathcal{A}u_i+\varepsilon$, where $\mathcal{A} : L^2(\mathbb{R}^2) \rightarrow L^2(\mathbb{R}^2)$ is the forward Fourier transform and by $\mathcal{A}^*:L^2(\mathbb{R}^2)\rightarrow L^2(\mathbb{R}^2)$ the backward inverse Fourier transform. 
$\mathcal{A}$ is a linear and continuous operator according to Plancherel's theorem. In this work, we assume to have real-valued images, therefore, our MRI forward operator takes only the real part of complex-valued measurements.

Let $\Omega$ be the image domain, a connected bounded open subset of $\mathbb{R}^2$ of class $\mathcal{C}^1$, and $u:\Omega \rightarrow \mathbb{R}$ be the sought single reconstructed image depicting the true underlying anatomy. For theoretical purposes, we assume that it is equal to $0$ on the boundary of $\Omega$, so that, we can extend it by $0$ and apply the operator $\mathcal{A}$.

We introduce  the unknown deformation as $\varphi_i : \bar{\Omega} \rightarrow \mathbb{R}^2$, between the $i$-th acquisition and the image $u$. Moreover, let $z_i$ be the associated displacements such that $\varphi_i=\mbox{Id}+z_i$, where $\mbox{Id}$ is the identity function. At the practical level, these deformations should be with values in $\bar{\Omega}$. Yet, at the theoretical level, working with such spaces of functions result in losing the structure of vector space. Nonetheless, we can show that our model retrieves deformations with values  in $\bar{\Omega}$ -  based on Ball's results~\cite{corona-bib:ball}. A deformation is a smooth mapping that is topology-preserving and injective, except possibly on the boundary if self-contact is allowed. We consider  $\nabla \varphi_i : \Omega \rightarrow M_2(\mathbb{R})$ to be the gradient of the deformation, where $M_2(\mathbb{R})$ is the set of real square matrices of order two. After introducing the set of notations for this work, we now turn to to motivate the introduction of our mathematical model.

\subsection{Joint Variational Model}

In this work, we seek to extract - \textit{jointly} -  from a set of MR acquisitions $x_i$ which are corrupted by motion, a clean static free motion-corrected reconstructed image $u$, along with the physiological dynamics through the deformation maps $\varphi_i$. 
With this purpose in mind, we consider the motion correction issue as a registration problem. In a unified variational framework, $u$ and all $\varphi_i$'s for $i=1,\cdots,T$, are seen as optimal solutions of a specifically designed cost function composed of: (1) a regularisation on the deformations $\varphi_i$ which prescribes the nature of the deformation, (2) a regularisation on $u$ that allows for removing artefacts whilst keeping fine details and  (3) a fidelity term measuring the discrepancy between the deformed mean - to which the forward operator has been applied -  and the acquisitions, intertwining then the reconstruction and the registration tasks.

In order to allow for large and smooth deformations, we propose viewing the shapes to be matched as isotropic, homogeneous  and hyperelastic materials. This outlook dictates the regularisation on the deformations $\varphi_i$, which is based on the stored energy function of an Ogden material. We introduce the following notations : $A : B =tr(A^TB)$ is the matrix inner product whilst $\|A\|=\sqrt{A:A}$ is the associated Frobenius matrix norm. In two dimensions, the stored energy function of an Ogden material, in its general form, is given by the following expression. $W_O(F)=\underset{i=1}{\overset{M}{\sum}}a_i\|F\|^{\gamma_i}+\Gamma(\mbox{det}F)$, with $a_i>0$, $\gamma_i\geq 1$ for all $i=1,\cdots,M$ and $\Gamma\, :\, ]0;\infty[ \rightarrow \mathbb{R}$, is a convex function satisfying $\underset{\delta \rightarrow 0^+}{\lim} \Gamma(\delta)=\underset{\delta \rightarrow +\infty}{\lim}\Gamma(\delta)=+\infty$. 

In this work, we focus on the particular energy, which reads: 

$W_{Op}(F)=\left\{ \begin{array}{cc}a_1\|F\|^4 + a_2\left( \mbox{det}F - \frac{1}{\mbox{det}F}\right)^4 &\text{ if } \mbox{det}F>0\\+\infty&\text{ otherwise},\end{array} \right.$ 

\noindent
with $a_1>0$, and $a_2>0$. The first term penalises the changes in length, whilst the second term enforces small changes in area. We check that this function falls within the general formulation of the stored energy function of an Ogden material: $W_{Op}(F) = \tilde W(\xi,\delta)=\left\{ \begin{array}{cc}a_1\|\xi\|^4 + a_2\left( \delta - \frac{1}{\delta}\right)^4 &\text{ if } \delta>0,\\+\infty&\text{ otherwise},\end{array} \right.$ 
$\tilde W : M_2(\mathbb{R}) \times \mathbb{R} \rightarrow \mathbb{R}$ is continuous since $\underset{\delta \rightarrow 0^+}{\lim}\tilde W(\xi,\delta)=\underset{\delta \rightarrow +\infty}{\lim}\tilde W(\xi,\delta)=+\infty$, and is convex with respect to $\delta$ \Big($g: \left| \begin{array}{c}\mathbb{R}^+_* \rightarrow \mathbb{R} \\ x \mapsto \left(x-\frac{1}{x}\right)^4 \end{array} \right.$ and $g'' : \left| \begin{array}{c} \mathbb{R}^+_* \rightarrow \mathbb{R} \\ x\mapsto 4\left(x-\frac{1}{x}\right)^2\left(\frac{3x^4+4x^2+5}{x^4}\right)\geq0\end{array}\right.$\Big). The design of the function $\Gamma$ is driven by Ball's results~\cite{corona-bib:ball} guaranteeing the deformation to be a bi-H\"older homeomorphism, and therefore, preserving the topology. It also controls that the Jacobian determinant remains close to one to avoid expansions or contractions that are too large.

The aforementioned regulariser is then applied along a discrepancy measure, which allows to interconnect the reconstruction and the registration tasks, and a regularisation of our free motion-corrected slice based on the weighted total variation.
With this purpose, let $g : \mathbb{R}^+ \rightarrow \mathbb{R}^+$ be an edge detector function satisfying $g(0)=1$, $g$ strictly decreasing and $\underset{r\rightarrow +\infty}{\lim}g(r)=0$. For the sake of simplicity, we  set $g_i = g( \|\nabla G_\sigma *\mathcal{A}^*x_i\|)$, with $G_\sigma$ being a Gaussian kernel of variance $\sigma$ and for theoretical purposes, we assume that there exists $c>0$ such that $0<c\leq g_i\leq 1$, and $g_i$ is Lipschitz continuous for each $i=1,\cdots,T$. In practice, we choose $g$ to be the Canny edge detector. 
We follow Baldi's arguments (\cite{debroux-bib:baldi}) to introduce the weighted $BV$-space and the associated weighted total variation related to the weight $g_i$, for each $i=1,\cdots,T$.
\begin{definition}[{\cite[Definition 2]{debroux-bib:baldi}}]
Let $w$ be a weight function satisfying some properties (defined in \cite{debroux-bib:baldi} and fulfilled here by $g_i$). We denote by $BV_w(\Omega)$ the set of functions $u\in L^1(\Omega,w)$, which are integrable with respect to the measure $w(x)dx$, such that:
\begin{align*}
 \sup\left\{ \int_\Omega u \mbox{div}(\varphi)\,dx\,:\, |\varphi|\leq w \text{ everywhere },\, \varphi\in Lip_0(\Omega,\mathbb{R}^2) \right\}<+\infty,
\end{align*}
with $Lip_0(\Omega,\mathbb{R}^2)$ the space of Lipschitz continuous functions with compact support. We denote by $TV_w$ the previous quantity.
\end{definition}

Given $E$ a bounded open set in $\mathbb{R}^2$, with boundary of class $\mathcal{C}^2$, then we have $TV_{g_i}(\xi_E)=|\partial E|(\Omega,g_i)=\int_{\Omega \cap \partial E} g_i\,dH^1$, where $\xi_E$ is the characteristic function of the set $E$. This quantity can be viewed as a new definition of the curve length- with a metric depending on the observations $x_i$. Minimising it is equivalent to locating the curve on the edges of $\mathcal{A}^*x_i$ where $g_i$ is close to $0$. We thus introduce the following fidelity term and regulariser for our single reconstruction:
\begin{align*}
 F(u,(\varphi_i)_{i=1,\cdots,T})= \frac{1}{T} \underset{i=1}{\overset{T}{\sum}}\delta TV_{g_i}(u \circ \varphi_i^{-1})+\frac{1}{2}\|\mathcal{A}((u\circ \varphi_i^{-1})_e)-x_i\|_{L^2(\mathbb{R}^2)}^2,
\end{align*}
where $\varphi_i^{-1}$ is the inverse deformation 
and $(u\circ \varphi_i^{-1})_e$ is the extension by $0$ of $(u\circ \varphi_i^{-1})$ outside $\Omega$. Indeed, we assume that $u=0$ and $\varphi_i=\varphi_i^{-1}=\mbox{Id}$ so that $u\circ \varphi_i^{-1}=0$ on the boundary $\partial \Omega$. Moreover, if $u\circ \varphi_i^{-1} \in BV_{g_i,0}(\Omega)$, with $BV_{g_i,0}(\Omega)=BV_0(\Omega) \cap BV_{g_i}(\Omega)$ (see \cite[Chapter 6.3.3]{corona-bib:demengel} for more details on the trace operator on the space of functions of bounded variations), then $u\circ \varphi_i^{-1}$ is null on the boundary $\partial \Omega$ and can be extended by $0$ outside the domain $\Omega$.
Due to the embedding theorem (\cite[Chapter 6.3.3]{corona-bib:demengel}), we have that $(u\circ\varphi_i^{-1})_e\in L^2(\mathbb{R}^2)$ and then Plancherel's theorem gives us $\mathcal{A}((u\circ \varphi_i^{-1})_e)\in L^2(\mathbb{R}^2)$ leading to the well-definedness of the fidelity term.
The first term of $F$ aims at aligning the edges of the deformed $u$ ($u\circ \varphi_i^{-1}$) with the ones of the different acquisitions, whilst regularising the reconstructed image $u$. The second quantity targets to make the Fourier transform of the deformed $u$ ($u\circ \varphi_i^{-1}$) as close as possible to the acquisitions $x_i$ (in the $L^2$ sense) - inspired by the classical fidelity term in MRI reconstruction models.

By combining both terms, we finally get the following minimisation problem
\begin{align}
 \inf I(u,(\varphi_i)_{i=1,\cdots,T}) = F(u,(\varphi_i)_{i=1,\cdots,T})+\frac{1}{T}\underset{i=1}{\overset{T}{\sum}}\int_\Omega W_{Op}(\nabla \varphi_i)\,dx. \tag{P} \label{initial_problem}
\end{align}

We now describe the numerical resolution of our joint model (\ref{initial_problem}).
\subsection{Numerical Method: Optimisation Scheme}
\label{sec:opt}
In order to tackle the nonlinearity of the problem, we introduce three auxiliary variables $v_i,\, w_i,f_i$ to mimic $\nabla \varphi_i$, $u\circ \varphi_i^{-1}$ and $w_i$, respectively, through quadratic penalty terms. This corresponds to reformulating the equality constraints $v_i=\nabla \varphi_i$, $w_i=u\circ \varphi_i^{-1}$ and $w_i=f_i$, into an unconstrained minimisation problem.
We then propose solving the following decoupled (discretised) problem
\begin{equation}
\begin{aligned}
    \min_{u,\varphi_i,v_i,w_i,f_i}\,& \frac{1}{T}\sum_{i=1}^T \underset{x \in \Omega}{\sum} W_{Op}(v_i(x)) + \frac{\gamma_1}{2} \|v_i - \nabla \varphi_i\|_2^2 \\
    &+ \frac{\gamma_3}{2} \|\mathcal{A}w_i -x_i \|_2^2 + \frac{\gamma_2}{2} \|(w_i - u \circ \varphi_i^{-1}) \sqrt{\operatorname{det} \nabla \varphi_i^{-1}}\|_2^2 \\
    &+ \frac{1}{2\theta}\|f_i -w_i\|_2^2 +\operatorname{TV}_{g_i}(f_i).
    \end{aligned}
    \label{eq:jointModel}
\end{equation}
We can now minimise \eqref{eq:jointModel} with an alternating splitting approach, where we solve the subproblem according to each individual variable whilst keeping the remaining variables fixed. %
We iterate the following updates for $k=0,1,\dots$. \smallskip

\textbf{Sub-problem 1: Optimisation over $v_i$.}  In practice, $v_i$ simulates the gradient of the displacements $z_i=(z_{i,1},z_{i,2})$ associated to the deformations $\varphi_i$. For every $v_i$, we have $v_i=\begin{pmatrix} v_{11} ~ v_{12} \\ v_{21} ~ v_{22}\end{pmatrix}$.  For the sake of readability, we drop here the dependency on $i$. We solve the Euler-Lagrange equation with a semi-implicit finite difference scheme (for $N$ iterations  and all $i$) and update $v_i$ as: 
\begin{equation*}
\begin{cases}
\begin{aligned}
v_{11}^{k+1}&=\frac{1}{1+dt \gamma} \Bigg( v_{11}^{k} + dt(-4a_1 \|v_i^{k}\|^2_F(v_{11}^{k}+1) - 4 a_2 (1+v_{22}^{k})  c_0 c_1
+ \gamma \frac{ \partial z_1^{k}}{\partial x} \Bigg), \\
v_{12}^{k+1}&=\frac{1}{1+dt \gamma} \Bigg( v_{12}^{k} + dt(-4a_1 \|v_i^{k}\|^2_F v_{12}^{k} + 4  a_2 v_{21}^{k} c_0 c_1
+ \gamma \frac{ \partial z_1^{k}}{\partial y} \Bigg), \\
v_{21}^{k+1}&=\frac{1}{1+dt \gamma} \Bigg( v_{21}^{k} + dt(-4a_1 \|v_i^{k}\|^2_F v_{21}^{k} + 4  a_2 v_{12}^{k} c_0 c_1
+ \gamma \frac{ \partial z_2^{k}}{\partial x} \Bigg), \\
v_{22}^{k+1}&=\frac{1}{1+dt \gamma} \Bigg( v_{22}^{k} + dt(-4a_1 \|v_i^{k}\|^2_F(v_{22}^{k}+1) - 4 a_2 (1+v_{11}^{k})  c_0 c_1+ \gamma \frac{ \partial z_2^{k}}{\partial y} \Bigg), \\
    \end{aligned}
    \end{cases}
\end{equation*}
with $c_0=\bigg(\operatorname{det}v_i^{k}- \frac{1}{\operatorname{det}v_i^{k}}\bigg)^3$ and $c_1=1+ \frac{1}{(\operatorname{det}v_i^{k})^2}$.

\smallskip
\textbf{Sub-problem 2: Optimisation over $\varphi_i$.} We solve the Euler-Lagrange equation in $\varphi_i$, after making the change of variable $y=\varphi_i^{-1}(x)$ in the $L^2$ penalty term, for all $i$, using an $L^2$ gradient 
flow scheme with a semi-implicit Euler time stepping (for N iterations in the same inner loop as $v_i$):
\begin{equation*}
    0=\gamma_1 \Delta \varphi_i^{k+1} + \gamma_1 \begin{pmatrix}\operatorname{div}v_{i,1}^{k+1} \\ \operatorname{div}v_{i,2}^{k+1} \end{pmatrix} + \gamma_2 (w_i^k \circ \varphi_i^{k} - u^k) \nabla w_i^k(\varphi_i^{k}).
\end{equation*}

\smallskip
\textbf{Sub-problem 3: Optimisation over $w_i$.} The update in $w_i$, for all $i$, reads:
\begin{equation*}
\begin{aligned}
w_i^{k+1} = \mathcal{A^{*}} \Bigg\{ \frac{\mathcal{A}\big(\gamma_2 \operatorname{det}\nabla (\varphi_i^{-1})^{k+1} u \circ (\varphi_i^{-1})^{k+1} + \frac{f_i^k}{\theta} \big)+ \gamma_3 x_i}{\gamma_2\operatorname{det}\nabla(\varphi_i^{-1})^{k+1}+\gamma_3+\frac{1}{\theta}} \Bigg\}.
\end{aligned}
\end{equation*}

\smallskip
\textbf{Sub-problem 4: Optimisation over $f_i$.} For each  $k$ and $i$, this is solved via Chambolle projection algorithm \cite{Chambolle2004}, as in \cite{Bresson2007}, for an inner loop over $n=0,1,\dots$
\begin{equation*}
    \begin{aligned}
    f_i^{n+1}&=w_i^{k+1}-\theta \operatorname{div}p_i^n, \\
    p_i^{n+1}&=\frac{p_i^n+ \delta_t \nabla (\operatorname{div} p_i^n - w_i^{k+1}/\theta)}{1+\frac{\delta_t}{g_i} | \nabla(\operatorname{div} p_i^n - w_i^{k+1}/\theta)| }.
    \end{aligned}
\end{equation*}
After enough iterations, we set $f_i^{k+1}=f_i^{n+1}$.

\smallskip
\textbf{Sub-problem 5: Optimisation over $u$.} Finally, using the same change of variables as in sub-problem 2, the problem in $u$ is simply
\begin{equation*}
    u^{k+1}=\frac{1}{T}\sum_{i=1}^T w_i^{k+1} \circ \varphi_i^{k+1}.
\end{equation*}

In the following section, we present the results obtained with our algorithmic approach, where the registration problem in $v_i$ and $\varphi_i$, $i=1,\dots,T$, has been solved in a multi-scale framework, from coarser to finer grids. It is worth noticing that - even though in $W_{Op}$, we have a term controlling that the Jacobian determinant remains positive, we introduce an additional step in our implementation, to monitor that the determinant does not become negative. 
To do so, we use the following regridding technique:
if at iteration $N$, $\operatorname{det}(\nabla \varphi^{N})$ becomes too small or negative, we save  $\varphi^{k,N-1}$ and re-initialise $\varphi^{k,N}=\operatorname{Id}$ and $v^{k,N}=0$, and $w^k = w^k \circ \varphi^{k,N-1}$, and we then continue the loop on $N$. Finally, $\varphi^{k,final}$ is the composition of all the saved deformations.
\setcounter{footnote}{0}
\section{Numerical Experiments}\label{sec:numerical_experiments} 
This sections describes the set of experimental results that are performed to validate our proposed joint model. \smallskip

\textbf{Data Description.} We evaluate our framework on two datasets. The \textit{Dataset I}\footnote{\url{http://www.vision.ee.ethz.ch/~organmot/chapter_download.shtml}} is a 4DMRI data acquired during free-breathing of the right liver lobe~\cite{Siebenthal2007}. It was acquired on a 1.5T Philips Achieva system, TR = 3.1 ms, coils =4, slices =25,  matrix size = $195\times166$. \textit{Dataset II}\footnote{\url{https://zenodo.org/record/55345\#.XBOkvi2cbUZ}} is a 2D T1-weighted dataset~\cite{BAUMGARTNER201783} acquired during free breathing of the entire thorax. It was acquired with a  3T Philips Achieva system  with matrix size = $215 \times 173$,  slice thickness= 8mm, TR=3.1ms and TE=1.9ms. The experiments reported in this section were run under the same condition in a CPU-based Matlab implementation. We used an Intel core i7 with 4GHz and a 16GB RAM.  \smallskip

\begin{wraptable}{l}{0.5\textwidth}
\vspace{-15pt}
    \resizebox{0.5\textwidth}{!}{  
    \begin{tabular}{c|c|c|c|c|c|c|c|c|c|c|c|}
    \cline{2-11}
    &   \cellcolor[HTML]{EFEFEF} $a_1$ & \cellcolor[HTML]{EFEFEF}$a_2$ & \cellcolor[HTML]{EFEFEF}$\gamma_1$ & \cellcolor[HTML]{EFEFEF}$\gamma_2$ & \cellcolor[HTML]{EFEFEF}$\gamma_3$ & \cellcolor[HTML]{EFEFEF}$\theta$ & \cellcolor[HTML]{EFEFEF}$\sigma$ & \cellcolor[HTML]{EFEFEF}$k$ & \cellcolor[HTML]{EFEFEF}$N$ & \cellcolor[HTML]{EFEFEF}$n$\\
    \hline
    \cellcolor[HTML]{EFEFEF}Dataset I     & 1& 50& 5& $10^5$ &30 &5 &2& 2& 500&500  \\
    \hline
    \cellcolor[HTML]{EFEFEF}Dataset II      & 1& 50 &5 & $10^5$ &15 & 5 &1.5 &2 & 500& 500 \\
    \hline
    \end{tabular}
    }
    \caption{Parameter values.}
    \label{tab:my_label}
\vspace{-20pt}
\end{wraptable} 
\textbf{Parameter Selection.} 
We present a discussion on the role of parameters in our model. The ranges of these parameters are rather stable for both experiments as shown in Table~\ref{tab:my_label}.
The parameters in the regularisation term for the registration $W_{Op}$ are $a_1$ and $a_2$. The former controls the smoothness of the deformation whilst the latter is more sensitive and can be viewed as a measure of rigidity allowing for a trade-off between topology preservation and the capacity of handling large deformations. 
If $a_2$ is chosen too small, the deformation may lose its injectivity property resulting in the loss of topology preservation, while if chosen too big, the registration accuracy may be impaired. 
The parameters $\gamma_1$,$\gamma_2$ and $\theta$, balance the $L^2$ penalty terms. 
Finally, $\gamma_3$ weights the data-fitting term for the reconstruction. For the FAIR approach - which is described next - we set the regularisation parameter for Dataset I$ =0.1$ and DatasetII$ =1$. 

\begin{figure*}[t!]
\centering
\includegraphics[width=1\textwidth]{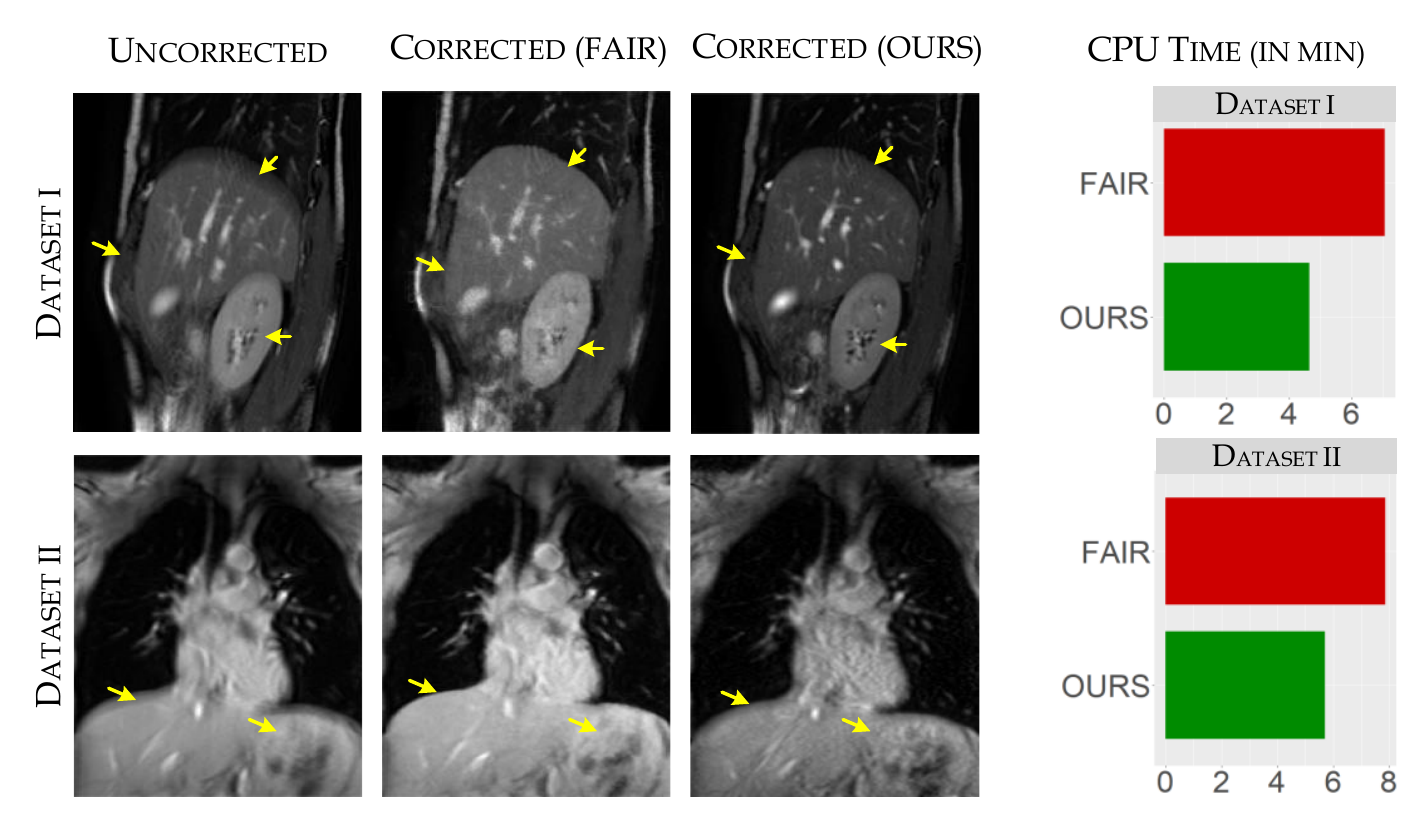}
\caption{Free-motion corrected comparisons: (From left to right) coarse Euclidean mean reconstructed image, FAIR and our joint model. Visual assessment in terms of blurring artefacts, loss of details and contrast are pointed out with yellow arrows. Elapsed time, in min, comparison between FAIR and ours approaches. }
\label{fig::resFig1}
\end{figure*}
\begin{figure*}[t!]
\centering
\includegraphics[width=1\textwidth]{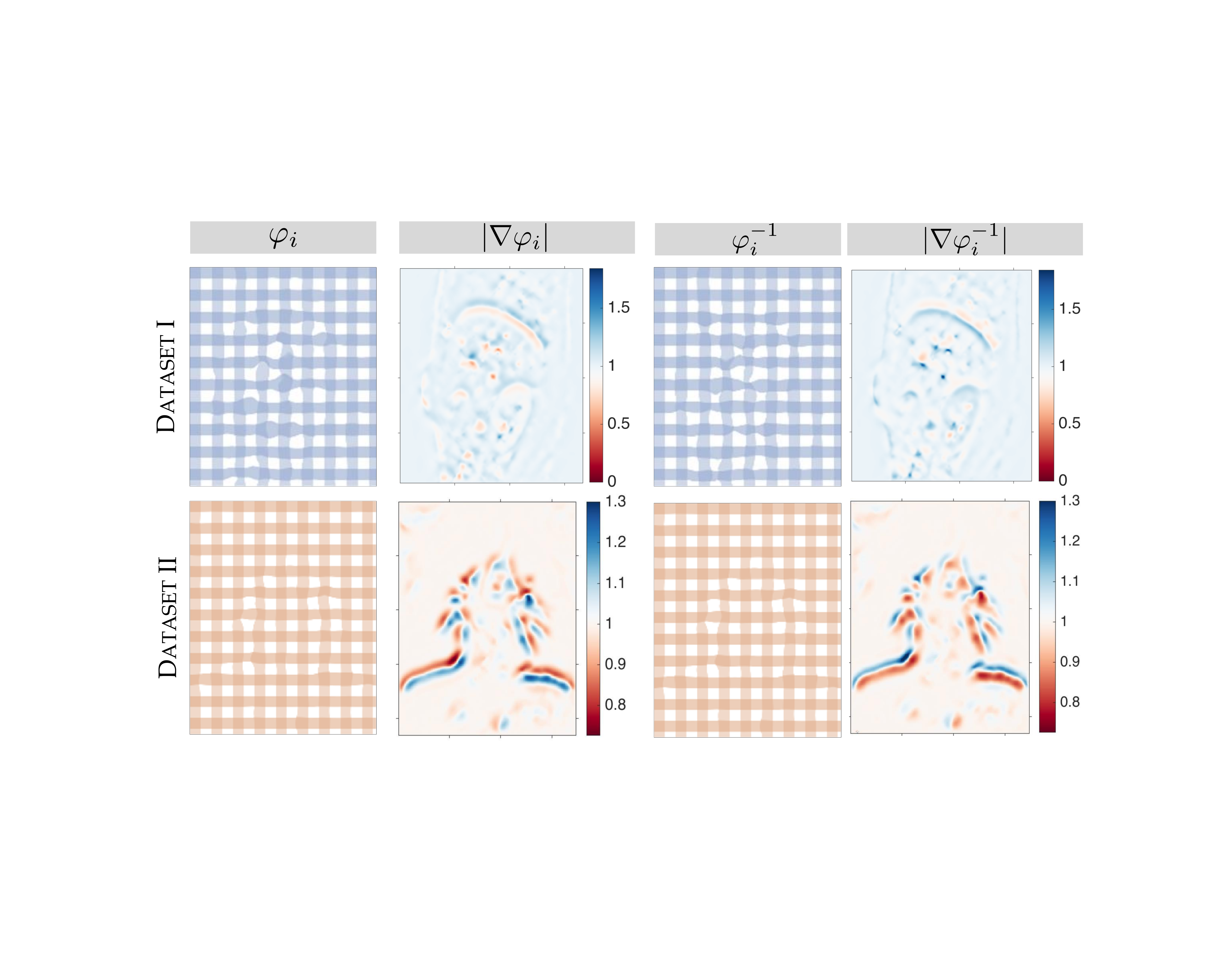}
\caption{Estimated motion and determinant maps of the deformation Jacobian. This is shown for the transformation $\varphi_i$ and its inverse $\varphi_i^{-1}$ for our two datasets.}
\label{fig::resFig2}
\end{figure*}
\smallskip
\textbf{Comparison with FAIR.} We begin by evaluating our joint model against a sequential approach using the FAIR toolbox~\cite{2009-FAIR}, which is a well-established framework for image registration. For a fair comparison, we set the characteristics of the compared approach closer to ours, as follows. Firstly, we choose the sum of squared differences as a distance measure on  $A^*x_i$, and secondly, we set the hyperelastic regulariser - as it allows for similar modeling of the deformations.

The computations of FAIR and our approach were performed under the same conditions.
We ran the same multi-scale approach. Moreover, we set the same reference image that we selected for the initialisation of our alternating scheme. 
This is important because of the non-convexity of our joint problem. After registering each frame to the reference, we average the results to obtain a mean image. In Fig.~\ref{fig::resFig1}, we show the outputs - for both datasets - \textit{uncorrected} being the coarse Euclidean mean, FAIR referring to the results of the process described above; and \textit{ours} being the reconstructed image $u$ obtained with our model. 
By visual inspection, we observe that the coarse Euclidean mean displays highly blurry effects along with loss of details whilst our approach improves in this regard. 
We also observed that - although FAIR improves over the coarse mean - outputs generated with FAIR have a remarkable loss of contrast  and detail. In addition, we point out that even though our method needs to solve 5 subproblems, the final computational time is less than the sequential FAIR approach.

\begin{figure*}[t!]
\centering
\includegraphics[width=1\textwidth]{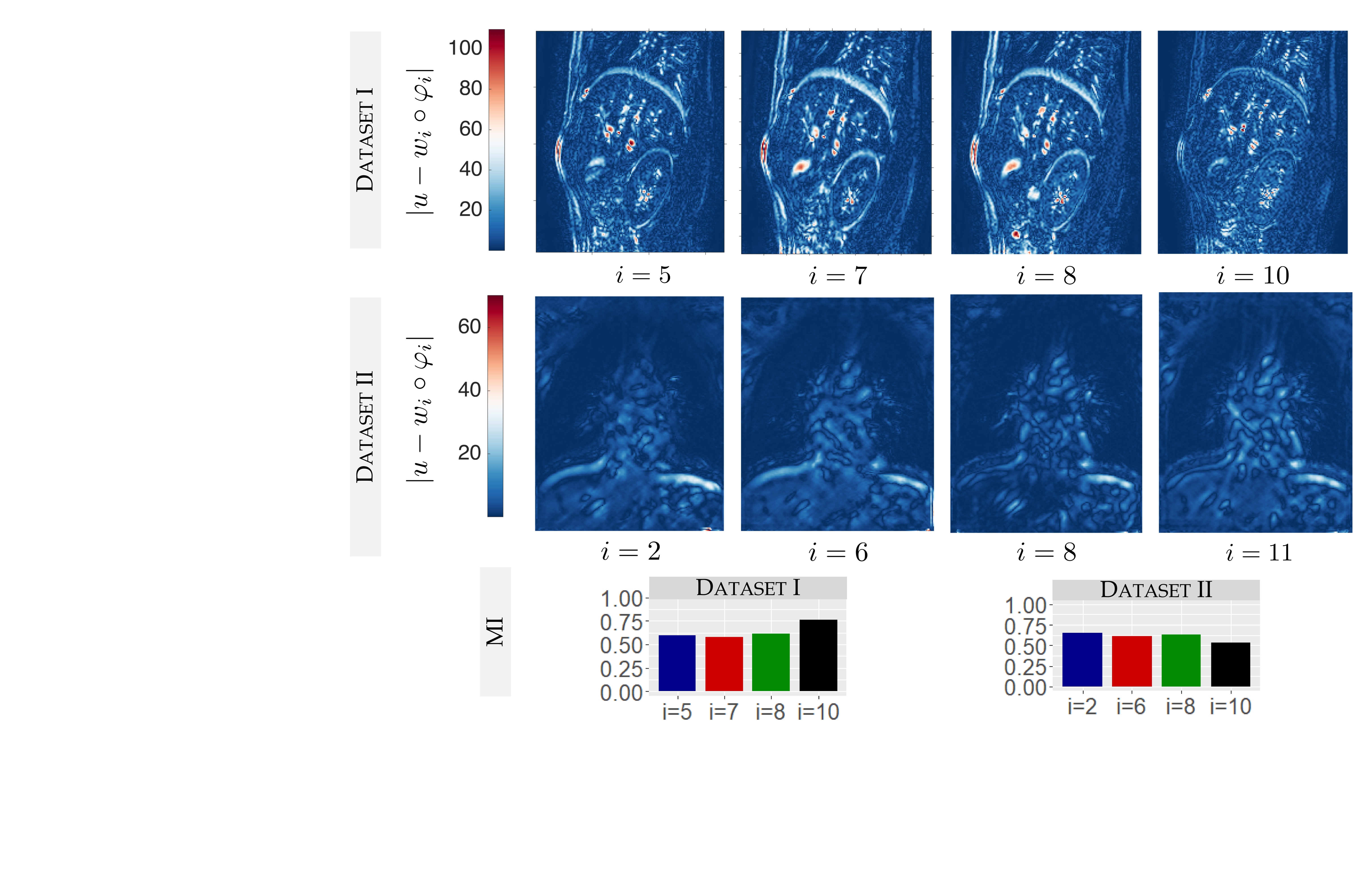}
\caption{ Top and middle rows: visual assessment of the difference map, at specific times,  between the resulting reconstruction $u$ (in the range $[0,255]$) and registered raw data $w_i \circ \varphi_i$. Bottom row: quantitative evaluation using MI values between $u$ and $w_i \circ \varphi_i$ for the same frames. }
\label{fig::resFig3}
\end{figure*}
Overall, visual inspection of the results shows the following drawbacks in which our joint model improves: (1) blurry effects: for both datasets, our model gives clearer images compared to the other approaches in which small details and shapes are lost, and (2) contrast preservation: our method is able to better preserve the contrast than other approaches, especially the one based on FAIR in which loss of contrast is clearly visible.

\textbf{Further Analysis of our Approach.} In Fig.~\ref{fig::resFig2}, we can observe that our method produces a reasonable estimation of the breathing dynamics, where the motion fields are visualised by a deformation grid. More precisely, it displays - for the two datasets - the estimated  motion $\varphi_i$ and its inverse $\varphi_i^{-1}$ for the frame $i=5$ along with their corresponding Jacobian determinant maps, whose values are interpreted as follows: small deformations when values are closer to $1$,  big expansions when values are greater than $1$, and big contractions when values are smaller than $1$.
Moreover, one can observe that the determinants remain positive in all cases, that is to say, our estimated deformations are physically meaningful and preserve the topology as required in a registration framework.

In Fig.~\ref{fig::resFig3}, we show further results for the registration.
We display a difference map (top) between our (adaptive) reference $u$, and the registered data $w_i \circ \varphi_i$ for sample frames for Dataset I and II. 
Moreover, we show the mutual information (MI) between the same images as a quantitative measure of performance (bottom), 
where larger MI values indicates better alignment. Since the $L^2$ fidelity term drives the registration process in our optimisation framework, it sounds relevant to use the difference map to illustrate the accuracy of our model. However, one should also be aware that a small change in contrast - which often happens between two frames and thus between $u $ and $w_i \circ \varphi_i$ - is clearly visible on the difference map, even if the overall anatomies are well aligned. Keeping in mind that, for the sake of tractability, the parameters are set globally and not adjusted for each frame, we can see that the intensities and the structures are well aligned for all the examples. This is corroborated by large MI values.

\section{Conclusion}
In this work we presented a new mathematical model to tackle motion correction in MRI in a joint reconstruction/registration framework. The underlying idea is to exploit redundancy in the temporal resolution in the data to compensate for motion artefacts due to breathing. The proposed multi-tasking approach solves simultaneously the problems of reconstruction and registration to achieve higher image quality and accurate estimation of the breathing dynamics. We mathematically motivated our joint model and presented our optimisation scheme. Furthermore, our experimental results show the potential of our approach for clinical applications. 

\section*{Acknowledgements} VC acknowledges the financial support of the Cambridge Cancer Centre. Support from the CMIH and CCIMI University of Cambridge is greatly acknowledged.

\bibliographystyle{splncs03}
\bibliography{refs}

\end{document}